\documentclass{elsart}
\usepackage{amssymb}
\usepackage{amsmath}
\usepackage[english]{babel}
\usepackage{graphicx}
\usepackage{theorem}
\newtheorem{theorem}{Theorem}

\newtheorem{corollary}[theorem]{Corollary}

\newtheorem{proposition}[theorem]{Proposition}
\newtheorem{remark}[theorem]{Remark}
\begin{document}
\begin{frontmatter}
\title{Equilibrium and eigenfunctions estimates in the semi-classical regime.}
\author{Brice Camus}
\thanks[label2]{Work supported by the
\textit{SFB/TR12} \textit{Symmetries\&Universality in Mesoscopic
Systems}.}
\address{Ruhr-Universit\"at Bochum, Fakult\"at f\"ur Mathematik,\newline%
Universit\"atsstr. 150, D-44780 Bochum, Germany.\newline Email :
brice.camus@univ-reims.fr}
\begin{abstract}
\noindent  We establish eigenfunctions estimates, in the
semi-classical regime, for critical energy levels associated to an
isolated singularity. For Schr\"odinger operators, the asymptotic
repartition of eigenvectors is the same as in the regular case,
excepted in dimension 1 where a concentration at the critical
point occurs. This principle extends to pseudo-differential
operators and the limit measure is the Liouville measure as long
as the singularity remains integrable.
\end{abstract}
\end{frontmatter}

\section{Introduction.}
The problem we consider here concerns the asymptotic behavior of
eigenvectors of a self-adjoint operator and follows the works of
Colin de Verdi\`ere \cite{CdV} and Zelditch \cite{Zel}, on the
basis of a result stated first by Schnirelman. We are more
precisely interested in a proof of a microlocal concentration
phenomena, in the semi-classical regime, as established in
\cite{BPU,CdV-P}. The adaptation to semi-classical analysis was
done in \cite{HMR}, following a technic proposed by Voros
\cite{Vor}. We also mention \cite{D-H} for a more general approach
in the scattering setting.

Consider a quantum operator $P_h$, realized as a self-adjoint
operator acting on a dense subset of $L^2(\mathbb{R}^n)$. A
typical example, studied in section 2, is the Schr\"odinger
operator $P_h=-h^2\Delta+V$ where the potential $V$ is smooth and
bounded from below. If the spectrum of $P_h$ is discrete in
$[E-\varepsilon ,E+\varepsilon ]$, a sufficient condition for this
is given below, we can enumerate the eigenvalues in this interval
as a sequence $\lambda _{j}(h)$ with finite multiplicities. We
note $\psi_j^h$ the corresponding normalized eigenvectors, i.e.
\begin{equation*}
P_h \psi_j^h=\lambda_j(h) \psi_j^h,\text{ } ||\psi_j^h||_{L^2}=1.
\end{equation*}
Our objective is to establish eigenfunctions estimates:
\begin{equation*}
\left\{
\begin{matrix}
\nu_j(a)=\left\langle \mathrm{Op}^w_h(a)\psi_j^h,\psi_j^h \right\rangle,\\
\lambda_j(h)\rightarrow E, \text{ }h\rightarrow 0^+,
\end{matrix}
\right.
\end{equation*}
where $a\in S^{0}(\mathbb{R}^{2n})$, so that by the
Calderon-Vaillancourt Theorem:
\begin{equation*}
f\mapsto \mathrm{Op}^w_h(a)f(x)=\frac{1}{(2\pi
h)^n}\int\limits_{\mathbb{R}^{2n}}
a(\frac{x+y}{2},\xi)e^{\frac{i}{h} \langle x-y,\xi\rangle}
f(y)dyd\xi,
\end{equation*}
is bounded on $L^2(\mathbb{R}^n)$. Note that the statement of the
problem is local w.r.t. $E$ and we are here interested in the case
of $E=E_c$ critical. Each $\nu_j(a)$ measures the observable
$\mathrm{Op}_h^w(a)$ in the state $\psi_j^h$. Let
$\Phi_t=\exp(tH_p)$ be the Hamiltonian flow of the principal
symbol $p$ of $P_h$. Interpreted as distributions, these measures
are $\Phi_t$-invariant which easily follows from:
\begin{equation*}
\left\langle
e^{-\frac{it}{h}P_h}\mathrm{Op}^w_h(a)e^{\frac{it}{h}P_h}\psi_j^h,\psi_j^h
\right\rangle=\left\langle
\mathrm{Op}^w_h(a)e^{\frac{it}{h}\lambda_j(h)}
\psi_j^h,e^{\frac{it}{h}\lambda_j(h)}\psi_j^h
\right\rangle=\nu_j(a).
\end{equation*}
By Egorov's Theorem
$e^{-\frac{it}{h}P_h}\mathrm{Op}^w_h(a)e^{\frac{it}{h}P_h}$ is an
operator of principal symbol $(a\circ \Phi_t)$ and $\nu_j$ is
invariant under $\Phi_t$, up to $\mathcal{O}(h)$.

We recall that $E$ is regular if $\nabla p\neq 0$ on the energy
surface:
\begin{equation*}
\Sigma _{E}=\{(x,\xi )\in T^{*}\mathbb{R}^{n}\text{ }/\text{
}p(x,\xi )=E\},
\end{equation*}
and critical otherwise. When $E_c$ is critical, $\Sigma_{E_c}$ is
not a smooth manifold. For $E$ regular, $\Sigma_E$ inherits a
measure, invariant by $\Phi_t$, given by:
\begin{equation*}
\mathrm{dLvol}(z)=\frac{dz}{||\nabla p(z)||} \text{ }_{|\Sigma_E},
\text{ } z\in \Sigma_E,
\end{equation*}
where $dz$ is the Riemannian surface element. We note
$\mathcal{V}(E)$ the associated volume of $\Sigma_E$ and we obtain
a probability measure via:
\begin{equation*}
d\mu^{E}(z) = \frac{1}{\mathcal{V}(E)}\mathrm{dLvol}(z).
\end{equation*}
Note that if $E_c$ is critical, $d\mu^{E_c}$ has a sense if and
only if $1 \in L^1 (\Sigma_{E_c},\mathrm{dLvol})$.

Via a wave equation approach, substituting here the heat equation
strategy of \cite{CdV} on a compact manifold, the problem is
related at the first order to the geometry of the energy surfaces.
Accordingly, the integrability of $\mathrm{dLvol}$ has a strong
effect on the asymptotic behavior of the sequence $\nu_j$.

The micro-local concentration near a singularity was studied in
\cite{CdV-P} in dimension one for a non-degenerate instable
equilibrium where they proved that $\nu_j$ converges to the Dirac
mass at the equilibrium. We are first interested in the case of
Schr\"odinger operators, but a generalization to
pseudo-differential operators provides more examples. Finally, we
recall that in Riemannian geometry, e.g. for compact surfaces of
negative curvature, the question to know if the full sequence
converges to the invariant measure (quantum unique
ergodicity) is still an open problem.\medskip\\
\textbf{Definitions.}\\
We define now the objects used in sections 2\&3. If $E_c$ is a
critical energy level, we pick an $h$-dependant interval:
\begin{equation}
I(h)=[E_c-dh,E_c+dh],\text{ } d>0.
\end{equation}
The associated counting function of eigenvalues is:
\begin{equation}
\Upsilon(h)= \sum\limits_{\lambda_j(h)\in I(h)} \left\langle
\psi_j^h, \psi_j^h \right\rangle =\# \{j \text { / }
\lambda_j(h)\in I(h) \}.
\end{equation}
For $A=\mathrm{Op}_h^w(a)$ a pseudodifferential operator of order
zero, whose principal symbol is $a\in S^0(\mathbb{R}^{2n})$ we
put:
\begin{equation}
\Upsilon_a(h)= \sum\limits_{\lambda_j(h)\in I(h)} \left\langle A
\psi_j^h, \psi_j^h \right\rangle=\sum\limits_{\lambda_j(h)\in
I(h)} \nu_j(a).
\end{equation}
Observe that $\Upsilon(h)=\Upsilon_1(h)$ since for every
quantization used in this contribution the symbol of the identity
is 1.
\section{Schr\"odinger operators.}%
Let $p(x,\xi)=\xi^2 +V(x)$, where the potential $V\in
C^{\infty}(\mathbb{R}^n)$ is real valued. To obtain a well defined
spectral problem, we use:
\begin{quote}
$(\mathcal{H}_{1})$ \textit{There exists } $C\in \mathbb{R}$ \textit{ such that: }%
$\liminf\limits_{\infty}V >C$.
\end{quote}
By a classical result, $P_h=-h^2\Delta+V(x)$ is essentially
self-adjoint. Note that $(\mathcal{H}_1)$ is always satisfied if
$V$ goes to infinity at infinity. Let $J=[E_1,E_2]$, with
$E_2<\liminf\limits_{\infty}V$. Since $p^{-1}(J)$ is compact, the
spectrum $\sigma (P_{h})\cap J$ is discrete and consists in a
sequence $\lambda _{1}(h)\leq \lambda _{2}(h)\leq ...\leq \lambda
_{j}(h)$ of eigenvalues of finite multiplicities, if $h$ is small
enough. Next, we impose the singularity:
\begin{quote}
$(\mathcal{H}_{2})$\textit{ On }$\Sigma _{E_c}$ \textit{the symbol
}$p$\textit{ has an isolated critical point }$z_{0}=(x_{0},0).$
\textit{This critical points can be degenerate but is associated
to a local extremum of $V$
\begin{equation*}
V(x)=E_c+ V_{2k}(x)+\mathcal{O}(||x-x^j_0||^{2k+1}), \text{
}k\in\mathbb{N}^{*},
\end{equation*}
where $V_{2k}$, homogeneous of degree $2k$, is definite positive
or negative.}
\end{quote}
The case $k=1$, i.e. a non-degenerate singularity in dimension
$n$, is treated in \cite{BPU} without any extremum condition. In
Fig.1 the line in bold is the critical energy level attached to
the top of a one dimensional symmetric degenerate double well.
Observe the unstability near the recurrent critical point.
\begin{figure}[h!]
\centering
\includegraphics{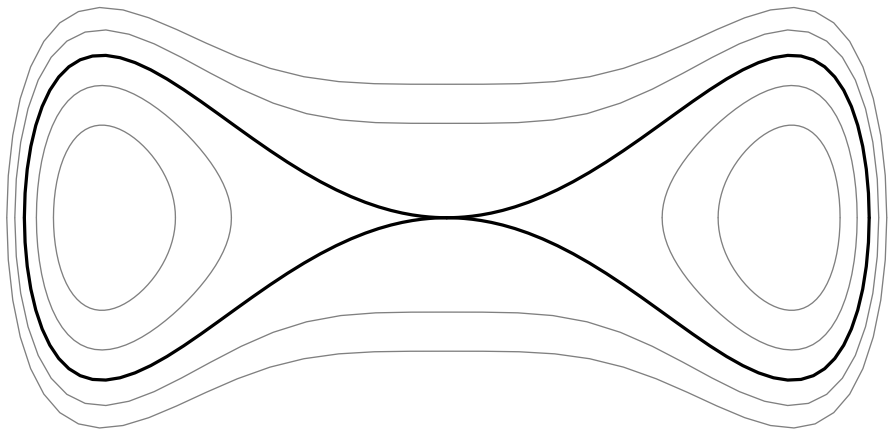}
\caption{Energies surfaces of $V(x)=-x^4+x^6$.}
\end{figure}

To simplify notations we write $z=(x,\xi)\in T^{*}\mathbb{R}^n$
and $z_0$ for a critical point. The first result concerns the
statistical behavior of the sequence $\nu_j(a)$.
\begin{theorem}\label{Main}
Assume $(\mathcal{H}_{1})$ and $(\mathcal{H}_{2})$ satisfied. If
$n>1$ we have
\begin{equation*}
\lim\limits_{h\rightarrow 0^+}\frac{ \Upsilon_a (h)}{\Upsilon(h)}
= \int a d\mu^{E_c},
\end{equation*}
but in dimension 1 we obtain:
\begin{equation*}
\lim\limits_{h\rightarrow 0^+} \frac{ \Upsilon_a (h)}{\Upsilon(h)}
=\langle \delta_{z_0},a\rangle =a(x_0,0).
\end{equation*}
\end{theorem}
These relations are statistical since when $h\rightarrow 0$ the
number of eigenvalues in $I(h)$ tends to infinity, see Proposition
\ref{estimates counting}. We define:
\begin{equation*}
K(h)=\{j\in\mathbb{N} \text{ / } \lambda_j(h)\in I(h)\}.
\end{equation*}
Following the results of \cite{BPU,CdV,Zel,HMR}, for $n>1$ Theorem
\ref{Main} implies that if $\Phi_t$ is ergodic on $\Sigma_{E_c}$
there exists a density one subset $L(h)\subset K(h)$ such that for
all integer valued function $h\rightarrow j(h)\in L(h)$ we have:
\begin{equation*}
\lim\limits_{h\rightarrow 0} \langle
\psi_{j(h)}^h,\mathrm{Op}^w_h(a)
\psi_{j(h)}^h\rangle=\int\limits_{\Sigma_{E_c}} a d\mu^{E_c}.
\end{equation*}
Here, density one simply means that:
\begin{equation*}
\lim\limits_{h\rightarrow 0} \, \frac{ \# L(h)} {\# K(h)}=1.
\end{equation*}
Since there is nothing new to prove, we refer to \cite{BPU,HMR}
for a precise study. More interesting, is the generalization of
the result of \cite{CdV-P}:
\begin{corollary}
In dimension 1, assume that $\Sigma_{E_c}$ is connected. Then, for
the weak * topology, we have $\nu_{j(h)} \rightarrow \delta(z_0)$,
$j(h)\in K(h)$, uniformly as $h\rightarrow 0$.
\end{corollary}
Hence, we obtain a concentration at $z_0$. Observe that
$\delta(z_0)$ is $\Phi_t$-invariant. If $a$ is simply a function,
the quantum probability satisfies:
\begin{equation*}
\lim\limits_{h\rightarrow 0}\, \lambda_j(h)=E_c \Rightarrow
\lim\limits_{h\rightarrow 0}\, |\psi_{j(h)}^h|^2(x)=\delta_{x_0}.
\end{equation*}
The interpretation is as follows. In dimension one, the invariant
measure on $\Sigma_{E_c}$ has a singularity:
\begin{equation*}
d\mathrm{Lvol}(z)\sim c \frac{dz}{||z-z_0||},\text{ near }
z_0=(x_0,0).
\end{equation*}
The measure is not integrable and the result has to be different.
For $n>1$, the singularity is integrable and an isolated critical
point has no effect. This reinforce the universality of the
Liouville measure in quantum ergodicity. But the case $n=1$ is
important since many problems, with symmetries, can be reduced to
the study of such a singular Schr\"odinger equation.
See e.g. \cite{CdV-P} for an application in Riemannian-geometry.\medskip\\
\textit{Preliminary remarks.}\\
The case of a local minimum of $V$ is not really deep. Since
$\xi^2\geq 0$, $z_0$ is a local extremum of $p$ and is an isolated
point of $\Sigma_{E_c}$. According to the results of \cite{Cam3},
the contribution of a minimum is significative only if $n=1$. We
consider now the non-trivial case of a local maximum of $V$,
corresponding to an unstable equilibrium of the flow. Finally,
since we use below the functional calculous, $p$ has to be a
symbol. But with $(\mathcal{H}_1)$ we can eventually modify the
potential $V$ outside of a compact, without modifying the main
results. Hence no extra assumption is required. Similar comments apply for section 3.\medskip\\
\textbf{Proof of Theorem \ref{Main}.}\\
We use the semi-classical trace formula technic. This approach,
analogous to the trace of the heat operator of \cite{CdV}, uses
the propagator $e^{\frac{it}{h}P_h}$ and a generalization of the
Poisson summation formula for this operator. Let $\varphi\in
\mathcal{S}(\mathbb{R})$, a Schwartz function. To approximate
$\Upsilon$ we define:
\begin{equation}\label{Def trace} \gamma (E_c,h,\varphi )=\sum\limits_{|\lambda
_{j}(h)-E_c| \leq \varepsilon}\varphi (\frac{\lambda
_{j}(h)-E_c}{h}).
\end{equation}
This object can be treated by mean of Fourier integral operators
(F.I.O.), see e.g. \cite{D-Hor}, and we refer to
\cite{BPU,Cam1,Cam3} for a detailed study of the trace. We recall
that the Tauberian approximation concerns expressions:
\begin{equation*}
\Upsilon^{\alpha}_{E,h}(\varphi)=\sum\limits_{j} \alpha_j(h)
\varphi(\frac{\lambda_j(h)-E}{h}).
\end{equation*}
Under our assumptions, the behavior of $\Upsilon^{\alpha}_{E,h}$
determines the behavior of the weighted counting functions:
\begin{equation*}
N_{E,d}^\alpha (h)= \sum\limits_{ |\lambda_j(h)-E|\leq dh}
\alpha_j(h).
\end{equation*}
Precisely, we will apply the results of section 6 of \cite{BPU}
for $\alpha_j(h)=1$ or $\nu_j(a)$. Strictly speaking,
$\alpha_j(h)\geq 0$ is required but, by a standard result of
pseudodifferential calculus, we can modify the quantization to
have $\nu_j(a)\geq 0$. This does not change the main results, see
Eq.(\ref{change quantization}) below.

To attain our objective, we can suppose that
$\mathrm{supp}(\hat{\varphi})\subset [-M,M]$, $M>0$. Let $\Theta
\in C_{0}^{\infty }(]E_{c}-\varepsilon ,E_{c}+\varepsilon \lbrack
)$, such that $\Theta =1$ in a neighborhood of $E_{c}$ and $0\leq
\Theta \leq 1$ on $\mathbb{R}$. We localize the problem near
$E_{c}$ by writing:
\begin{equation*}
\gamma (E_{c},h,\varphi)=\gamma _{1}(E_{c},h,\varphi)+\gamma
_{2}(E_{c},h,\varphi),
\end{equation*}
with:
\begin{gather*}
\gamma _{2}(E_{c},h,\varphi)=\sum\limits_{|\lambda _{j}(h)-E_c|
\leq \varepsilon}\Theta (\lambda _{j}(h))\varphi (\frac{\lambda
_{j}(h)-E_{c}}{h})\\
=\mathrm{Tr}\, \Theta(P_h) \varphi (\frac{P_h-E_{c}}{h}).
\end{gather*}
Where the last equality holds by support considerations. By a
classical result, see e.g. \cite{Cam1} Lemma 1, the term
$\gamma_1=\gamma-\gamma_2$ satisfies :
\begin{equation}
\gamma _{1}(E_{c},h,\varphi)=\mathcal{O}(h^{\infty }), \text{ as }
h\rightarrow 0.\label{S1(h)=Tr}
\end{equation}
The Fourier inversion formula for $\gamma_2$ and the previous
estimate provide:
\begin{equation}
\gamma (E_{c},h,\varphi)=\frac{1}{2\pi }\mathrm{Tr}\int\limits_{\mathbb{R}}e^{i%
\frac{tE_{c}}{h}}\hat{\varphi}(t)\mathrm{exp}(-\frac{i}{h}%
tP_{h})\Theta (P_{h})dt+\mathcal{O}(h^\infty).
\end{equation}
Next, with a function $\Psi \in C_{0}^{\infty }(T^{\ast
}\mathbb{R}^{n})$, with $\Psi =1\text{ near }z_{0}$, we write:
\begin{equation*}
\gamma_{2}(E_{c},h,\varphi)=\gamma_{z_0}(E_{c},h,\varphi)+\gamma
_{\mathrm{reg}}(E_{c},h,\varphi),
\end{equation*}
where:
\begin{equation}\label{local problem}
\gamma_{z_0}(E_{c},h,\varphi) =\frac{1}{2\pi }\mathrm{Tr}\int\limits_{\mathbb{R}}e^{i%
\frac{tE_{c}}{h}}\hat{\varphi}(t)\Psi ^{w}(x,hD_{x})\mathrm{exp}(-\frac{i}{h}%
tP_{h})\Theta (P_{h})dt,
\end{equation}
and $\gamma_{\mathrm{reg}}$ is simply the difference. The
micro-local term $\gamma_{z_0}$ contains the contribution of the
singularity and the discussion below determines if this term is
dominant. For finitely many critical point on $\Sigma_{E_c}$, we
could repeat the procedure. For the convenience of the reader, we
recall the contributions of an equilibrium to the trace formula.
\begin{proposition}\label{trace Schrodinger}
Assume $(\mathcal{H}_1)$, $(\mathcal{H}_2)$ and that
$\hat{\varphi}\in C_{0}^{\infty}([-M,M])$, $M\leq M_0$. If $x_0$
is a local maximum of $V$ we have:
\begin{equation*}
\gamma _{z_{0}}(E_{c},h,\varphi)\sim
h^{-n+\frac{n}{2}+\frac{n}{2k}}
\sum\limits_{m=0,1}\sum\limits_{j,l\in\mathbb{N}^2}
h^{\frac{j}{2}+\frac{l}{2k}}\mathrm{log}(h)^m \Lambda
_{j,l,m}(\varphi ).
\end{equation*}
If $\frac{n(k+1)}{2k}\in \mathbb{N}$ and $n$ is odd then the
top-order term is :
\begin{equation*}
C_{n,k} \log (h)h^{-n+\frac{n}{2}
+\frac{n}{2k}}\int\limits_{\mathbb{S}^{n-1}}
|V_{2k}(\eta)|^{-\frac{n}{2k}} d\eta  \int\limits_{\mathbb{R}}
|t|^{n\frac{k+1}{2k}-1} \varphi (t) dt.
\end{equation*}
Otherwise the first non-zero coefficient are given by:
\begin{equation*}
h^{-n+\frac{n}{2} +\frac{n}{2k}}%
\left\langle T_{n,k},\varphi\right\rangle%
\int\limits_{\mathbb{S}^{n-1}}
|V_{2k}(\eta)|^{-\frac{n}{2k}}d\eta.
\end{equation*}
\end{proposition}
This result is the contribution of an equilibrium since the
distributional coefficients have a non-discrete support, contrary
to the Weyl-term supported in $t=0$ and the contributions of
closed orbits supported by the length spectrum. The distributions
$T_{n,k}$ and the universal constants $C_{n,k}\neq 0$ depend only
on $(n,k)$ and are explicitly determined in \cite{Cam4}. Mainly,
we need the order w.r.t. $h$ of these contributions determined by
the functions:
\begin{equation*}
w_c(h)=h^{-n+\frac{n}{2}+\frac{n}{2k}}\log (h)^j,\text{ }j=0
\text{ or } 1.
\end{equation*}
First, we give a natural application of Proposition \ref{trace
Schrodinger}.
\begin{proposition}\label{estimates counting}
The microlocal counting function satisfies:
\begin{equation*}
\Upsilon(h)=%
\left\{
\begin{matrix}
2d \mathcal{V}(E_c)(2\pi h)^{1-n}+o(h^{1-n}),\text{ if } n>1,\\
\Lambda(\chi_{[-d,d]}) w_c(h)+o(w_c(h)), \text{ if } n=1.
\end{matrix}
\right.
\end{equation*}
Here $\Lambda$ is the first non-zero distribution of Prop.
\ref{trace Schrodinger} and $\chi_{[-d,d]}$ the characteristic
function of $[-d,d,]$.
\end{proposition}
\textit{Proof.} By construction we have:
\begin{equation*}
\gamma_{\mathrm{reg}}(E_{c},h,\varphi)=\frac{1}{2\pi }\mathrm{Tr}\int\limits_{\mathbb{R}}%
e^{i\frac{tE_{c}}{h}}\hat{\varphi}(t)(1-\Psi^{w}(x,hD_{x}))%
\mathrm{exp}(-\frac{i}{h}tP_{h})\Theta (P_{h})dt.
\end{equation*}
By the standard calculus on F.I.O. and an easy application of the
stationary phase method, as $h$ tends to 0 we obtain:
\begin{equation*}
\gamma_{\mathrm{reg}}(E_{c},h,\varphi)\sim\frac{\hat{\varphi}(0)}{(2\pi
h)^{n-1}} \mathrm{Lvol}(\Sigma_{E_c}\cap \mathrm{supp}(1-\Psi))
+\mathcal{O}(h^{2-n}).
\end{equation*}
Hence $\gamma_{\mathrm{reg}}$ always contributes at the order
$h^{1-n}$.\medskip\\
\textbf{Case of $n>1$}. We have $w_c(h)=o(h^{1-n})$ and :
\begin{equation*}
\mathrm{Lvol}(\Sigma_{E_c}\cap \mathrm{supp}(1-\Psi))\leq
\mathcal{V}(E_c)<\infty, \text{ } \forall \Psi.
\end{equation*}
It follows easily by shrinking the support of the cut-off $\Psi$
that:
\begin{equation}
\gamma(E_c,\varphi,h)=(2\pi h)^{1-n} \hat{\varphi}(0)
\mathcal{V}(E_c) +o(h^{1-n}).
\end{equation}
Since the distributional factor is :
\begin{equation*}
\hat{\varphi}(0)=\int\limits \varphi(t)dt=\left\langle 1, \varphi
\right\rangle,
\end{equation*}
replacing $\varphi$ by $\chi_{[-d,d]}$, via Theorem 6.3 of
\cite{BPU}, provides:
\begin{equation*}
\Upsilon(h)=2d \mathcal{V}(E_c)(2\pi h)^{1-n}+o(h^{1-n}).
\end{equation*}
\textbf{Case of $n=1$.} Here the contribution of the critical
point has a bigger order than the regular one. We obtain :
\begin{equation}
\gamma(E_c,\varphi,h)=w_c(h) \Lambda(\varphi) +o(w_c(h)).
\end{equation}
To apply the Tauberian argument of \cite{BPU}, we observe that our
distribution $\Lambda\in \mathcal{S}'(\mathbb{R})$ can be
represented by an element of $L^1_{\mathrm{loc}}(\mathbb{R})$ and
hence can be extended as a linear form on $C_0(\mathbb{R})\cap
L^\infty(\mathbb{R})$. We obtain:
\begin{equation}
\Upsilon(h)\sim  \Lambda (\chi_{[-d,d]})w_c(h).
\end{equation}
Which provides the desired result for $n=1$.$\hfill{\blacksquare}$\medskip\\
As a matter of illustration, for $n=k=1$, we have:
\begin{equation}
\log(h) |V_{2k}(x_0)|^{-\frac{1}{2}}
\int\limits_{-d}^{d}dt=%
\frac{2d\log (h)}{|V''(x_0)|^{\frac{1}{2}}}.
\end{equation}
Which is the result established for $\Upsilon(h)$ in
\cite{BPU,CdV-P}. Observe that for $n=1$, $k=1$ is the only case
where a logarithm occurs and $\Upsilon(h)$ is slowly increasing
w.r.t. $k$ since for all $k>1$ :
\begin{equation}
\Upsilon(h)\sim C
h^{\frac{1}{2k}-\frac{1}{2}}|V^{(2k)}(x_0)|^{-\frac{1}{2k}}.
\end{equation}
\textbf{Eigenfunctions estimates.}\\
We recall how to derive eigenfunctions estimates from the trace
formula. First, to insert an observable $A$ changes almost
nothing. If $\Pi$ is the spectral projector on $[E_c-\varepsilon,
E+\varepsilon]$, computing the trace in the basis $\psi_j^h$ and
by cyclicity:
\begin{equation*}
\mathrm{Tr} \left( \Pi A \varphi(\frac{P_h-E_c}{h})\right
)=\sum\limits_{|\lambda_j(h)-E_c|\leq \varepsilon} \left\langle
A\psi_j^{h},\psi_j^h \right\rangle
\varphi(\frac{\lambda_j(h)-E_c}{h}).
\end{equation*}
Since $A$ is a bounded operator, if $\varphi\in
\mathcal{S}(\mathbb{R})$ we can again smooth the problem via an
energy cut-off $\Theta(P_h)$, with an error of order
$\mathcal{O}(h^\infty)$. Hence we can insert
$A=\mathrm{Op}^w_h(a)$ in Eq.(\ref{local problem}) and the results
of Prop. \ref{trace Schrodinger} are the same after multiplication
by $a(z_0)$. Similarly, the regular contribution changes via :
\begin{equation*}
\frac{\hat{\varphi}(0)}{(2\pi h)^{1-n}} \int\limits_{\Sigma_{E_c}}
a(z) (1-\Psi(z)) \mathrm{dLvol}(z).
\end{equation*}
By evaluation of the trace, we have :
\begin{equation*}
\sum\limits_{|\lambda_j(h)-E_c|\leq \varepsilon} \left\langle
A\psi_j^h, \psi_j^h\right\rangle
\varphi(\frac{\lambda_j(h)-E_c}{h}) \sim c_0(\varphi) w(h)m(a)+
o(w(h)),
\end{equation*}
where $w(h)$ changes only if $n=1$. By Theorem 6.3 of \cite{BPU}
we obtain :
\begin{equation*}
\sum\limits_{\lambda_j(h)\in I(h)} \left\langle
A\psi_j^h,\psi_j^h\right\rangle= m(a)w(h)+o(w(h)).
\end{equation*}
In particular this implies that:
\begin{equation}
\lim\limits_{h\rightarrow 0^+} \frac{1}{\Upsilon(h)}
\sum\limits_{\lambda_j(h)\in I(h)} \left\langle A\psi_j^h,
\psi_j^h\right\rangle=\lim\limits_{h\rightarrow 0^+}
\frac{\Upsilon_a(h)}{\Upsilon(h)}=\frac{m(a)}{m(1)}.
\end{equation}
Substituting the correct expressions for these measures we obtain:
\begin{quote}
- for $n>1$ : $m$ is a constant multiple of the Liouville
measure.\\
- for $n=1$ : $m$ is a multiple of the delta-Dirac distribution in
$z_0$.
\end{quote}
\textbf{Extraction of a subsequence.}\\
We chose $a\geq 0$ and modify the quantization. Different choices
are possible: Friedrichs quantization $\mathrm{Op}^F$ as in
\cite{BPU} or the anti-Wick quantization $\mathrm{Op}^{AW}$ as in
\cite{HMR}. These quantization are positive, i.e.
\begin{equation*}
a\geq 0 \Rightarrow \left\langle f ,\mathrm{Op}^{AW}(a)f
\right\rangle\geq 0,\text{ } \forall f\in
C_0^{\infty}(\mathbb{R}^n).
\end{equation*}
Since $\mathrm{Op}_h^w(a)-\mathrm{Op}_h^{AW}(a)$ is of order -1,
we obtain :
\begin{equation}\label{change quantization}
\left\langle \psi_j^h,
(\mathrm{Op}_h^w(a)-\mathrm{Op}^{AW}_h(a))\psi_j^h\right\rangle=\mathcal{O}(h),
\end{equation}
and we can work with this positive operator. For $n>1$, under the
condition that $\Phi_t$ is ergodic on $\Sigma_{E_c}$, the
extraction of a convergent subsequence of density one is the same
as in \cite{BPU,CdV,HMR} to which we refer for a detailed proof.
For $n=1$, if $\Sigma_{E_c}$ is connected, there is only one
probability measure invariant by $\Phi_t$ and the full sequence
converges to $\delta_{z_0}$. Once the result is established for a
positive symbol it can be extended by linearity to any $a\in
S^0(\mathbb{R}^{2n})$.
\section{Pseudo-differential operators.}
The case of pseudo-differential operators provides more explicit
examples. Let $P_{h}=\mathrm{Op}_{h}^{w}(p(x,\xi)$, obtained by
Weyl quantization, where the symbol $p$ is real valued and smooth
on $T^*\mathbb{R}^n$. In general position, one can also consider
$h$-dependent symbols $\sum h^{j}p_{j}(x,\xi )$, see \cite{HMR}.
But, to simplify, we consider only the homogeneous case. As above
we impose:
\begin{quote}
$(\mathcal{A}_{1})$\textit{There exists }$\varepsilon_{0}>0$ \textit{ such that }%
$p^{-1}([E_c-\varepsilon_{0},E_c+\varepsilon_{0}])$\textit{\ is
compact.}
\end{quote}
As in section 2, $\sigma (P_{h})\cap [E_{c}-\varepsilon
,E_{c}+\varepsilon ]$ is discrete. A fortiori, $(\mathcal{A}_1)$
insures that $\Sigma_{E_c}$ is compact. Next, we chose an homogeneous singularity :
\begin{quote}
$(\mathcal{A}_{2})$\textit{ On }$\Sigma _{E_c}$\textit{,
}$p$\textit{ has a unique critical point }$z_{0}=(x_{0},\xi
_{0})$\textit{ and near }$z_{0}$ :
\begin{equation*}
p(z)=E_{c}+\mathfrak{p}_{k}(z)+\mathcal{O}(||(z-z_{0})||^{k+1}),\text{
} k>2,
\end{equation*}
\textit{where }$\mathfrak{p}_{k}$\textit{\ is homogeneous of degree }$k$%
\textit{\ w.r.t. }$z-z_{0}$.
\end{quote}
Strictly speaking, one could consider $k=2$. But this case is
precisely treated in \cite{BPU}. The case of a critical point
which is not an extremum is technical because the singularity is
transferred on the blow up of $z_0$. To obtain a problem that can
be explicitly solved, we consider the following hypothesis
inspired by H\"ormander's real principal condition:
\begin{quote}
$(\mathcal{A}_{3})$ \textit{We have }$\nabla \mathfrak{p}_{k}\neq
0$\textit{ on the set }$C (\mathfrak{p}_k)=\{ \theta \in
\mathbb{S}^{2n-1}\text{ / } \mathfrak{p}_{k}(\theta )=0
\}$\textit{.}
\end{quote}
For example $\mathfrak{p}_3(x,\xi)=x^3-\xi^3$ is admissible and
$p(x,\xi)=x^3-\xi^3+x^4+\xi^4$ satisfies all our hypothesis for
$E_c=0$.
\begin{remark} \rm{With $(\mathcal{A}_{3})$, contrary to the case of a local extremum,
$z_0$ is not an isolated point of $\Sigma_{E_c}$ which imposes to
study the classical dynamic in a neighborhood of $z_0$.  The study
of singularities like in $(\mathcal{A}_{3})$ is detailed in
\cite{G-C} chapter 4 to which we refer concerning the
integrability of $\mathrm{dLvol}$.}
\end{remark}
As in section 2, it is sufficient to study the local problem
$\gamma _{z_{0}}$ defined in Eq.(\ref{local problem}). The
contributions to the trace formula are :
\begin{proposition} \label{trace pseudo}
Under $(\mathcal{A}_{1})$ to $(\mathcal{A}_{3})$, we have an
asymptotic expansion:
\begin{equation*}
\gamma_{z_0}(E_c,\varphi,h)\sim
h^{\frac{2n}{k}-n}\sum\limits_{m=0,1} \sum\limits_{j=0}^\infty
h^{\frac{j}{k}}\log(h)^m \Lambda_{j,m}(\varphi),
\end{equation*}
where the logarithms only occur when $(2n+j)/k\in \mathbb{N}^*$
and $\Lambda_{j,m}\in \mathcal{S}'(\mathbb{R})$.\\
As concerns the leading term we obtain:\\
(1) If $k>2n$ (non-integrable singularity on $\Sigma_{E_c}$) we
have :
\begin{equation*}
\gamma _{z_{0}}(E_{c},\varphi,h)\sim h^{\frac{2n}{k}-n}\Lambda
_{0,0}(\varphi )+\mathcal{O}(h
^{\frac{2n+1}{k}-n}\mathrm{log}(h)),\text{ as }h\rightarrow 0,
\end{equation*}
where $\Lambda _{0,0}$ is a universal distribution.\\
(2) If the ratio $2n/k\in\mathbb{N}$ we obtain logarithmic
contributions :
\begin{equation*}
\gamma _{z_{0}}(E_{c},\varphi,h)\sim
h^{\frac{2n}{k}-n}\mathrm{log}(h)\Lambda_{0,1} (\varphi
)+\mathcal{O}(h ^{\frac{2n}{k}-n}),\text{ as }h\rightarrow 0,
\end{equation*}
(3) For $2n>k$, $2n/k\notin \mathbb{N}$ the result is as in (1)
with a different distribution.
\end{proposition}
These results describe very precisely the singularity at $z_0$.
But this is not our purpose here and we refer to \cite{Cam5} for a
detailed formulation of these contributions. For $n=1$, $k=2$, the
case \textit{(2)} agrees with section 2 and allows to recover
some results established in \cite{BPU,CdV-P}.\medskip\\
\textbf{Application to microlocal measures.}\\
The proof is exactly the same as in section 2. The main difference
is that the singularity on $\Sigma_{E_c}$ can be of arbitrary
order. In our setting, according to Prop.\ref{trace pseudo} the
top order coefficient changes if and only if we have:
\begin{equation}
\frac{2n}{k}-n<1-n\Leftrightarrow \frac{2n}{k}<1.
\end{equation}
If $k<2n$ the singularity is integrable and contributes at a lower
order compared to $h^{1-n}\mathcal{V}(E_c)$. But if $k\geq 2n$,
which corresponds to a non-integrable singularity for
$d\mu^{E_c}$, the main term changes. To summarize, we obtain:
\begin{equation*}
\lim\limits_{h\rightarrow 0} \,\frac{ \Upsilon_a (h)}{\Upsilon(h)}
= \left\{
\begin{matrix}
\int a d\mu^{E_c}, \text{ for $k<2n$},\\
a(z_0),\text{ for $k\geq 2n$}.
\end{matrix}
\right.
\end{equation*}
Contrary to section 2, observe that for $k\geq 2n$ and if $n>1$ we
do not obtain the convergence of the full sequence $\nu_{j(h)}$,
$j(h)\in K(h)$, to the dirac-mass at the equilibrium. The
obstruction is that an invariant probability measure can be
supported by the closed orbits of $\Sigma_{E_c}$.\medskip\\
\textbf{Comments.}\\
From these 2 families of example the conclusion is that the
limiting measure changes only if $\Sigma_{E_c}$ carries a measure
such that $1\notin
L^1_{\mathrm{loc}}(\Sigma_{E_c},\mathrm{dLvol})$. Interpreted as a
quantum measurement, one obtain a very precise localization : if
$a=0$ around $z_0$ the limit is the Liouville-measure but if
$a(z_0)\neq 0$ the limit strongly differs.

An interesting problem would be to study the repartition in
presence of 2 equilibria $z_1$, $z_2$ on $\Sigma_{E_c}$ of the
same nature and with a non-integrable singularity. In this case
any convex combination :
\begin{equation}
\nu=a\delta(z_1)+(1-a)\delta(z_2),\text{ } a\in[0,1],
\end{equation}
provides an invariant probability measure. A natural question is
to determine if the limiting measures are equally distributed
between $z_1$ and $z_2$.

\end{document}